\documentclass[12pt]{article} 
\usepackage{amsmath,amssymb,eucal} 
\usepackage{bbm}  
\font\tenrm=cmr10

\font\cmssl=cmss10 at 12 pt  
\font\bigss=cmssdc10 scaled 2300

\font\cmsslll=cmss10 at 14 pt

\renewcommand{\b}{\beta}

\newcommand{\e}{\epsilon}  
  
\newcommand{\g}{\gamma}

\renewcommand{\k}{\kappa}  
\renewcommand{\l}{\lambda}  
\renewcommand{\o}{\omega}

\newcommand{\s}{\sigma}  
\renewcommand{\t}{\tau}

\newcommand{\G}{\Gamma}

\renewcommand{\S}{\Sigma}

\newcommand{\bR}{\mathbb{R}}  
\newcommand{\bZ}{\mathbb{Z}}

\newcommand{\bS}{\mathbb{S}}

\newcommand{\gh}{\mathfrak{h}}

\newcommand{\so}{\mathfrak{so}}  
\newcommand{\su}{\mathfrak{su}}  
\newcommand{\spin}{\mathfrak{spin}}  
\newcommand{\gsp}{\mathfrak{sp}}

\newcommand\Sp{\mathrm{Sp}}

\newcommand\SU{\mathrm{SU}}  
  
\newcommand\Spin{\mathrm{Spin}}  
\newcommand{\id}   {{\mathbbm{1}}}

\renewcommand{\square}{\kern1pt\vbox  
               {\hrule height 0.6pt\hbox{\vrule width 0.6pt\hskip 3pt  
    \vbox{\vskip 6pt}\hskip 3pt\vrule width 0.6pt}\hrule height0.6pt}  
    \kern1pt}

\newcommand{\ra}{\rightarrow}  
\DeclareMathOperator\tr{tr\,}  
  
\DeclareMathOperator\End{End\,}

\newtheorem{Th}{Theorem}  
\newtheorem{Prop}{Proposition}  
\newtheorem{Cor}{Corollary}  
\newtheorem{Lem}{Lemma}  
\newtheorem{Def}{Definition}  
  
\newcommand{\bt}{\begin{Th}\ \ }  
\newcommand{\et}{\end{Th}}  
\newcommand{\bp}{\begin{Prop}\ \ }  
\newcommand{\ep}{\end{Prop}}  
\newcommand{\bc}{\begin{Cor}\ \ }  
\newcommand{\ec}{\end{Cor}}  
\newcommand{\bl}{\begin{Lem}\ \ }  
\newcommand{\el}{\end{Lem}}  
\newcommand{\bd}{\begin{Def}\ \ }  
\newcommand{\ed}{\end{Def}}  
  
\newcommand{\pf}{\noindent{\it Proof:\ \ }}  
\newcommand{\qed}{\hfill\square}  
\newcommand{\n}{\nabla}  
  
\newcommand{\ot}{\otimes}

\newcommand{\be}{\begin{equation}}  
\newcommand{\ee}{\end{equation}}  
  
\newcommand\re[1]{(\ref{#1})}  
\newcommand{\arr}{\begin{array}{rlll}}  
\newcommand{\ea}{\end{array}}  
\newcommand{\bea}{\begin{eqnarray}}  
\newcommand{\eea}{\end{eqnarray}}  
\newcommand{\bean}{\begin{eqnarray*}}  
\newcommand{\eean}{\end{eqnarray*}}  
  
\catcode`@=11  
\@addtoreset{equation}{section}  
\catcode`@=12

\begin{document}  
\begin{titlepage}  
%\rightline{hep-th/yymmnnn}  
%\rightline{draft: \today}  
\rightline{} 
\vskip 1.5 true cm  
\begin{center}  
{\bigss  On pseudo-Riemannian manifolds\\[.5em] with many Killing spinors}  
\vskip 1.0 true cm   
{\cmsslll  D.V.\ Alekseevsky and  V.\ Cort\'es} \\[3pt] 
{\tenrm   The University of Edinburgh and Maxwell Institute for
Mathematical Sciences\\ 
JCMB, The King's buildings,  
Edinburgh, EH9 3JZ, UK \\
D.Aleksee@ed.ac.uk}\\[1em]  
{\tenrm   Department Mathematik 
und Zentrum f\"ur Mathematische Physik\\ 
Universit\"at Hamburg, 
Bundesstra{\ss}e 55, 
D-20146 Hamburg, Germany\\  
cortes@math.uni-hamburg.de}\\[1em]   
November 13, 2008 
\end{center}  
\vskip 1.0 true cm  
%%%%%%%%%%%%%%%%%%%%%%%%%%%%%%%%%%%%%%%%%%%%%%%%%%%%%%%%  
\baselineskip=18pt  
\begin{abstract}  
\noindent  
Let $M$ be a pseudo-Riemannian spin manifold of dimension 
$n$ and signature $s$ and denote by $N$ the rank of the real spinor bundle.
We prove that $M$ is locally homogeneous if it admits more 
than $\frac{3}{4}N$ independent Killing spinors with the same Killing number, 
unless  $n\equiv 1 \pmod 4$ and $s\equiv 3 \pmod 4$. We also prove that  
$M$ is locally homogeneous if it admits $k_+$ independent Killing spinors with 
Killing number $\l$ and $k_-$ independent Killing spinors with 
Killing number $-\l$ such that $k_++k_->\frac{3}{2}N$, 
unless $n\equiv s\equiv 3\pmod 4$.  
Similarly, a pseudo-Riemannian manifold 
with more than $\frac{3}{4}N$ independent      
\emph{conformal} Killing spinors 
is \emph{conformally} locally homogeneous. 
For (positive or negative) definite metrics, the bounds 
$\frac{3}{4}N$ and $\frac{3}{2}N$ in the above results 
can be relaxed to $\frac{1}{2}N$ and $N$, respectively. Furthermore, we prove 
that a pseudo-Riemannnian spin manifold with more than 
$\frac{3}{4}N$ parallel spinors is flat and that
$\frac{1}{4}N$ parallel spinors suffice if the metric is definite. 
Similarly, a Riemannnian spin manifold with more than $\frac{3}{8}N$ 
Killing spinors with the Killing number $\l \in \bR$ 
has constant curvature  $4\l^2$. For Lorentzian or negative definite
metrics the same is true with the bound $\frac{1}{2}N$.  
Finally,
we give a classification of (not necessarily complete) Riemannian 
manifolds admitting Killing spinors, which provides an inductive construction
of such manifolds.  
\end{abstract}

\end{titlepage}  
%\tableofcontents
\section*{Introduction}
Figueroa-O'Farrill, Meessen and Philip showed in \cite{FMP} that
M-theory backgrounds with more than 24 supersymmetries are locally 
homogeneous. Notice that 24 is $3/4$ of the maximal possible number of
independent supersymmetries, which is 32, the 
dimension of the spinor module of $\Spin(1,10)$. (Notice also that $11\equiv
3 \not\equiv 1\pmod{4}$.)   
This result is obtained from a careful analysis 
of the Killing spinor equations of M-theory. 

In this paper, inspired by the work of Figueroa-O'Farrill et al, 
we study Killing spinors in pseudo-Riemannian and 
conformal geometry  for arbitrary dimensions $n$ and signatures $s$.
We show that conformal Killing spinors give rise to conformal
Killing polyvectors and, under some simple assumptions,  
that Killing spinors give rise to Killing polyvectors, see 
Theorem \ref{1stThm}. More precisely, in equation \re{fundEqu}, 
we define  a $\wedge^kTM$-valued 
bilinear form 
\[ (s,t) \mapsto [s,t]_k, \]
on the spinor bundle of a pseudo-Riemannian spin manifold 
$(M,g)$, which to a pair of conformal Killing 
spinors $s,t$ associates a conformal Killing polyvector field 
$\o =[s,t]_k$. For $k=1$ we obtain conformal Killing vector 
fields. 

Using the above correspondence, 
we prove that the existence of more than $3/4$ of the maximal possible 
number $N$ of independent Killing spinors implies local homogeneity in the 
pseudo-Riemannian as well as in the conformal setting, see 
Theorem \ref{2ndThm} for the precise statement. 
For (positive or negative) definite metrics we prove that
more than $\frac{1}{2}N$ Killing spinors suffice to obtain the local 
homogeneity. 
In the pseudo-Riemannian (but not in the conformal) 
setting, our argument requires $n\not\equiv 1\pmod{4}$ or 
$s\not\equiv  3\pmod{4}$.  
Allowing imaginary ``Killing numbers'' $\l I\in \End S$, 
where $\l\in \bR$ and $I^2=-\id$, see \re{imagKSEqu}, we can 
prove a similar result also in the case $n\equiv 1\pmod{4}$, 
$s\equiv  3\pmod{8}$. In the remaining case, where $n\equiv 1\pmod{4}$ and 
$s\equiv  7\pmod{8}$, our method does not allow to obtain 
the local homogeneity from the existence of 
Killing spinors with the same Killing number. Instead we have to 
assume the existence of $k_+$ 
Killing spinors with Killing number $\l$ and $k_-$ Killing spinors
with Killing number $-\l$. If $k_++k_->\frac{3}{2}N$,  
then we prove that the 
pseudo-Riemannian manifold is locally homogeneous, provided that
$n\not\equiv 3\pmod{4}$ or $s\not\equiv 3\pmod{4}$. This covers,
in particular the case $n\equiv 1\pmod{4}$. For definite metrics the
assumption can be relaxed to $k_++k_->N$.  

Using the correspondence between Killing spinors on $(M,g)$ and
parallel spinors on the metric cone $(\hat{M},\hat{g})$ over $M$, 
see Definition \ref{coneDef} and Theorem \ref{BaerThm}, and our recent work 
\cite{ACGL} we are able to obtain more precise information for 
 Riemannian and Lorentzian manifolds. In fact, in  
Theorems \ref{parallelThm}, \ref{RiemThm} and \ref{Lor/ndThm}  we prove: 

\bt \label{introThm}
\begin{enumerate}
\item[(i)] 
A pseudo-Riemannian spin manifold with more than 
$\frac{3}{4}N$ linearly 
independent parallel spinors is flat. If the metric is definite,
then $\frac{1}{4}N$ parallel spinors suffice. 
\item[(ii)] A Riemannian spin manifold 
with more than $\frac{3}{8}N$
Killing spinors with the Killing number $\l \in \bR$ 
has constant curvature  $4\l^2$. 
\item[(iii)] A pseudo-Riemannian spin manifold with a negative definite or
Lorentzian metric with more than $\frac{1}{2}N$
Killing spinors with the Killing number $\l \in \bR$ 
has constant curvature  $4\l^2$. 
\end{enumerate}
\et 

Notice that a negative definite metric $g$ of positive 
scalar curvature $s$
corresponds to a positive definite metric $-g$ of negative scalar 
curvature $-s$.  
We also prove that a Riemannian spin manifold with $\frac{3}{8}N$ 
Killing spinors with the Killing number $\l \in \bR  \setminus \{ 0\}$ 
can be locally represented in the form
\[ M= I \times M_1 \times M_2,\quad g=ds^2 + \cos^2(s) g_1 + \sin^2(s) g_2,
\]
where $(M_1,g_1)$ is of constant curvature $1$ or of dimension $\le 1$,  
$(M_2,g_2)$ is a seven-dimensional $3$-Sasakian manifold and 
$I\subset (0,\frac{\pi}{2})$ is an intervall, 
see Theorem \ref{RiemThm}.   

In Theorem \ref{ClassThm}, 
we give a local 
classification of Riemannian 
manifolds admitting a nontrivial Killing spinor, which extends B\"ar's 
classification \cite{B} of Killing spinors on complete Riemannian manifolds. 

\section{From Killing spinors to Killing polyvectors}
Let $(M,g)$ be an $n$-dimensional pseudo-Riemannian manifold. We will always 
assume that $M$ is connected.  
\bd  A $k$-vector field  
$\o \in \G (\wedge^kTM)\cong \G (\wedge^kT^*M)$ ($k\ge 1$) is called {\cmssl Killing} if 
\[ X\lrcorner \n_X\o =0, \quad\mbox{for all}\quad X\in TM.\]
It is called {\cmssl conformally Killing} if there exists a $(k-1)$-vector field
$\tilde{\o}$ such that 
\begin{equation} X\lrcorner \n_X\o = g(X,X)\tilde{\o}, 
\quad\mbox{for all}\quad X\in TM.\label{confKvEqu}\end{equation}
\ed

\bp  \begin{enumerate}
\item[(i)] $\o\in \G (\wedge^kTM)$ is Killing if and only if 
$\dot{\g}\lrcorner \o$ is a parallel $(k-1)$-vector field along $\g$,  
for every geodesic $\g$: 
\begin{equation}
\n_{\dot{\g}}(\dot{\g}\lrcorner \o )=0. \label{conservEqu}
\end{equation}
\item[(ii)] Let $(M,g)$ be a pseudo-Riemannnian manifold with 
indefinite metric $g$. Then 
$\o$ is conformally Killing if and only if 
$\n_{\dot{\g}}(\dot{\g}\lrcorner \o )=0$,  
for every null geodesic $\g$.
\end{enumerate}
\ep 

\pf An obvious calulation shows that a (conformal) Killing polyvector $\o$
satisfies the equation \re{conservEqu} for every (null) geodesic $\g$. 
The converse statement in (i) is also clear, since every vector
$X$ is the velocity vector of a geodesic. 
To prove the converse statement in (ii), let  
$\eta \in \wedge^{k-1}T_pM$ and denote by $\beta$ the  
symmetric bilinear form such that  
$\eta \lrcorner (X\lrcorner \n_X\o) = \beta (X,X)$, for all $X\in T_pM$. 
By \re{conservEqu}, we have $\beta(X,X)=0$ for all $X$ in the null cone of 
$g$. This shows that $\beta$ is a multiple of $g_p$, since the null cone 
determines the indefinite scalar product $g_p$ up to scale,   
and implies \re{confKvEqu}.  
\qed

\noindent 
{\bf Remarks:} 1) 
For $k=1$ (i) reduces to the well know fact that the scalar product 
of a Killing vector field with the velocity vector of a geodesic  
is constant,
which was observed by Clairaut for surfaces of revolution. 
In virtue of (ii), 
conformal Killing polyvectors give rise to conservation laws 
in general relativity. 
In particular, the function $g(\dot{\g},Y)$ is constant
along any null geodesic $\g$ if $Y$ is a conformal Killing vector 
field.\\[.5em] 
2) It is easy to see that an $n$-vector field $\o$ on an $n$-dimensional manifold 
is conformally Killing if and only if it is parallel.\\[.5em]
3) The equation \re{confKvEqu} easily implies 
\[ \tilde{\o } = \frac{1}{n}\tr \n \o = 
\frac{1}{n}\sum g^{ij}e_i\lrcorner (\n_{e_j}\o),\]
where $e_i$ is any basis and 
$(g^{ij})$ is the matrix inverse to $g_{ij}=g(e_i,e_j)$.  

Let $(M,g)$ be a (strongly oriented) pseudo-Riemannian spin manifold and 
$S\ra M$ its (real) spinor bundle. 
\bd A spinor field $s\in \G(S)$ is called 
{\cmssl Killing} with Killing number $\lambda \in \bR$ if 
\[ \n_Xs=\lambda Xs,\quad\mbox{for all}\quad X\in TM,\]
where $Xs$ is the Clifford product of the vector $X$ and the
spinor $s$. 
It is called {\cmssl conformally Killing} if there exists a
spinor field $\tilde{s}\in \G(S)$ such that 
\begin{equation}  \n_Xs= X\tilde{s},\quad\mbox{for all}\quad X\in TM.
\label{confKsEqu}  
\end{equation}
\ed 

\noindent 
{\bf Remarks:} 1) Using the Clifford relation, $XY+YX=-2g(X,Y)$, 
the equation \re{confKsEqu} easily implies 
\begin{equation} \tilde{s} = -\frac{1}{n}Ds,\label{stildeEqu}
\end{equation} 
where $Ds=\sum g^{ij}e_i\n_{e_j}s$ is the Dirac operator. 
In particular, any Killing spinor is an eigenspinor for the 
Dirac operator: $Ds=-n\l s$.\\[.5em] 
2) The Killing number is related to the scalar curvature
by the formula $scal=4n(n-1)\l^2$. Therefore, the scalar curvature 
of a pseudo-Riemannian manifold which admits a Killing spinor 
is constant and the Killing numbers of different 
Killing spinors on the same manifold coincide up to a sign. 
It is well known that a Riemannian manifold which admits a Killing spinor
is Einstein, but this is no longer true for indefinite
pseudo-Riemannian manifolds, see [Bo] and references therein. 

We denote by $\g_v:S_p\ra S_p$ the Clifford multiplication with $v\in T_pM$
and define a linear map $\g : \wedge^kT_pM \ra \End (S_p)$, for all $k\ge 1$,  
by 
\[ \g_{v_1\wedge \cdots\wedge v_k} := \frac{1}{k!}\sum_{\s \in \mathfrak{S}_k}\e (\s) 
\g_{v_{\s 1}}\cdots
\g_{v_{\s k}},\]  
where $\mathfrak{S}_k$ is the symmetric group. 
For $\l \in \wedge^0T_pM=\bR$ we put $\g_\l = \l\id\in \End (S_p)$. 

A bilinear form $h$ on the spinor module satisfying 
\begin{eqnarray}
h (s,t) &=& \s h (t,s),\nonumber\\   
h (\g_Xs,t) &=& \t h (t,\g_Xs)\label{tauEqu}, 
\end{eqnarray} 
for all spinors 
$s,t$ and all vectors $X$,  
is called {\cmssl admissible} of {\cmssl symmetry} $\s$ and {\cmssl type} 
$\tau$, 
where $\s, \t \in \{-1,+1\}$. The admissible
bilinear forms on the spinor module were classified in \cite{AC}
and there always exists a \emph{nondegenerate} admissible
bilinear form. 
An admissible form is automatically invariant under the 
connected spin group and, hence, 
defines a parallel section of $S^*\ot S^*$. In the following,  
$h$ shall always denote a  parallel nondegenerate section of $S^*\ot S^*$ of 
symmetry $\s$ and type $\t$. Notice that \re{tauEqu}
implies 
\[ h(\g_\xi s,t)=\tau^k(-1)^{\frac{(k-1)k}{2}}h(s,\g_\xi t),\quad
\mbox{for all}\quad \xi \in \G (\wedge^kTM).\] 
Using the bilinear form $h$ we define, for $k\ge 1$, a parallel section
$[\cdot ,\cdot]_k\in \G (S^*\ot S^*\ot \wedge^kTM)$ by
\begin{equation} \label{fundEqu} g([s ,t]_k,\xi )=h(\g_\xi s,t)\quad 
\forall \xi \in  \G (\wedge^kTM), s,t\in \G (S).
\end{equation} 
(Here 
$g$ is canonically extended to a nondegenerate
symmetric bilinear form on the exterior algebra.)   
Such brackets occur in the classification of polyvector super-Poincar\'e 
algebras, see \cite{AC,ACDV}. For $k=0$ we put 
$[s,t]_0=h(s,t)$. 
\bt \label{1stThm} 
Let $s$, $t$ be conformal Killing spinors on an $n$-dimensional 
pseudo-Riemannian spin 
manifold $(M,g)$. Then $\o =[s,t]_k \in \G (\wedge^kTM)$ ($k\ge 1$) 
is a conformal 
Killing polyvector; 
\[ X\lrcorner \n_X\o =g(X,X)\tilde{\o}\quad \forall X\in TM,\]
where $\tilde{\o} \in \G (\wedge^{k-1}TM)$ is given by 
\begin{equation} 
n\tilde{\o} = (-1)^{k-1}[Ds,t]_{k-1}+\t [s,Dt]_{k-1}.\label{otildeEqu}
\end{equation}

\et 

\pf Let $(e_i)$ be a local frame and $\xi = X\wedge \eta$, where $X\in \G(TM)$, 
$\eta\in \G (\wedge^{k-1}TM)$ and 
$X\lrcorner\eta=0$. We shall assume that, at a given point $p\in M$,  
$\n X|_p =\n e_i|_p=0$ and $\n \eta|_p=0$. 
Then we compute at $p$:
\begin{eqnarray*}
g(\n_X\o , \xi ) &=&h(\g_\xi \n_Xs,t)+h(\g_\xi s,\n_X t)\\
&=& h(\g_\xi \g_X\tilde{s},t)+h(\g_\xi s,\g_X \tilde{t})\\
&=& -g(X,X)\left( (-1)^{k-1}h(\g_\eta\tilde{s},t)+ \t h(\g_\eta s,\tilde{t})
\right)\\
&=& -g(X,X)\left( (-1)^{k-1}g([\tilde{s},t]_{k-1},\eta) + \t 
g([s,\tilde{t}]_{k-1},\eta) \right) .
\end{eqnarray*} 
This implies that $\o$ is a conformal Killing polyvector and that
\begin{equation}\label{tildeEqu}
\tilde{\o} = (-1)^k[\tilde{s},t]_{k-1}-\t [s,\tilde{t}]_{k-1}.
\end{equation}  
Expressing $\tilde{s}, \tilde{t}$ by \re{stildeEqu}, we obtain \re{otildeEqu}. 
\qed 
\bc \label{mainCor} Let $s$ and $t$ be Killing spinors with Killing numbers  
$\lambda$ and $\mu$,  respectively, and $\o =[s,t]_k$. Then the following 
is true. 
\begin{enumerate}
\item[(i)] $\o$ is a conformal Killing polyvector
with $\tilde{\o} = (\l(-1)^k -\mu \t )[s,t]_{k-1}$.   
\item[(ii)] If $\mu=(-1)^k\t \l$, then $\o =[s,t]_k$ is a Killing polyvector. 
\item[(iii)] If $\l=\mu=0$ then $\o$ is parallel. 
\end{enumerate}
\ec
\section{Manifolds with many Killing spinors}
\bt \label{2ndThm} 
Let  $(M,g)$ be a pseudo-Riemannian spin manifold of dimension $n$,
signature $s$ and with spinor
bundle $S$ of rank $N$. 
\begin{enumerate}
\item[(i)] If $(M,g)$ admits $k>\frac{3}{4}N$ conformal
Killing spinors, which are linearly independent at $p\in M$, then  
$(M,g)$ admits $n$ conformal Killing vector fields, which are linearly 
independent at $p\in M$. 
\item[(ii)] Assume that $n\not\equiv 1 
\pmod 4$ or $s\not\equiv 3 \pmod 4$. If $(M,g)$ admits $k>\frac{3}{4}N$ 
Killing spinors with the same Killing number, 
which are linearly independent at $p\in M$, then  
$(M,g)$ admits $n$ Killing vector fields, which are linearly 
independent at $p\in M$.  
\item[(iii)] Assume that $n\equiv 1 \pmod 4$ and $s\equiv 3 \pmod 8$. 
Then $S$ admits a parallel hypercomplex structure $J_1,J_2,J_3=J_1J_2\in 
\G (\End S)$,  which commutes with Clifford multiplication. Let $I$ be any  
complex structure on $S$ which is a linear combination of $J_1,J_2,J_3$
with constant coefficients.  If $(M,g)$ admits $k>\frac{3}{4}N$ 
solutions $s\in \G (S)$ of the equation 
\begin{equation} \label{imagKSEqu}\n_Xs=\l XIs,\quad\mbox{for all}
\quad X\in TM,
\end{equation}
with the same $\l\in \bR$, which are linearly independent at $p\in M$, then 
$(M,g)$ admits $n$ Killing vector fields, which are linearly 
independent at $p\in M$.
\item[(iv)]   Assume that $n\not\equiv 3\pmod{4}$ or $s\not\equiv 3 \pmod 4$. 
If  $(M,g)$ admits $k_+$ 
Killing spinors with the Killing number $\l$, which are independent at $p$, 
and $k_-$ Killing spinors with the Killing number $-\l$, which are 
independent at $p$, such that $k_++k_->\frac{3}{2}N$, 
then it admits $n$ Killing vector fields, which are independent at $p$.  
\item[(v)] If $g$ is definite, then (i)-(iv) hold 
under the weaker assumptions 
$k>\frac{1}{2}N$ and $k_++k_->N$, respectively. 
\end{enumerate}
\et 

\pf $S$ carries a parallel nondegenerate 
bilinear form $h$ of symmetry $\s$ and type $\t$, 
see \re{tauEqu}. Moreover, there  
exists such a form of type $\tau=-1$, unless  $n\equiv 1 
\pmod 4$ and $s\equiv 3 \pmod 4$, see \cite{AC}. 
(The Pin$(n)$-invariant scalar product on the 
spinor module associated with a positive definite scalar product, for 
instance, has $\t=-1$.) 
By Theorem \ref{1stThm}, for any pair of conformal Killing spinors 
$s$, $t$, the vector field $[s,t]_1$ is conformal. Similarly, 
by Corollary \ref{mainCor}, if 
$s$, $t$ are Killing spinors with the same 
Killing number and $\t=-1$,   
then $[s,t]_1$ is a Killing vector field. 
Therefore, to prove (i) and (ii)  
it suffices to show that
\[ \Pi := [\cdot ,\cdot ]_1|_{S_0\ot S_0} : S_0\ot S_0 \ra T_pM\]  
is surjective if the subspace  $S_0\subset S_p$ 
spanned by the values of the given 
(conformal) Killing spinors at $p$ has dimension $>\frac{3}{4}\dim S_p$.
Suppose first that $g$ is definite.  Then we have to show
that $\Pi$ is surjective if
$\dim S_0 > \frac{1}{2}\dim S_p$.  By the definition of
$\Pi$, surjectivity is equivalent to: $\nexists v\in T_pM\setminus \{ 0\}$ 
such that $\g_vS_0\subset S_0^\perp$.  Suppose that there exists
$v\in T_pM\setminus \{ 0\}$ such that 
$\g_v|_{S_0} : S_0 \ra S_0^\perp$. If $\dim  S_0 > \frac{1}{2}\dim S_p$, then
$\dim S_0^\perp <\frac{1}{2}\dim S_p < \dim S_0$ and, thus, 
$\ker \g_v \neq 0$. Since $\g_v^2=-g(v,v)\id$, this implies 
$g(v,v)=0$ and, hence, $v=0$. This proves the surjectivity of 
$\Pi$, if $g$ is definite
and $\dim S_0 > \frac{1}{2}\dim S_p$. If $g$ is indefinite, we can 
only conclude that $v$ is a null vector. 
\bl \label{Lemma} 
For any non-zero null vector $v$ the subspace $L_v :=\ker \g_v = 
{\rm im}\, \g_v\subset S_p$ is 
$h$-isotropic of dimension $\frac{1}{2}\dim S_p$. 
\el 

\pf {}From $\g_v^2=0$ we get ${\rm im}\, \g_v \subset \ker \g_v$. 
Let $u$ be an other null vector such that $g(u,v)=1$. Then
${\rm im}\, \g_u \subset \ker \g_u$ and $\g_u\g_v+\g_v\g_u=-2\id$
implies $\ker \g_v \subset {\rm im}\, \g_v$ and, hence, 
$\ker \g_v = {\rm im}\, \g_v$. Therefore,
$\dim S_p-\dim \ker \g_v= \dim {\rm im}\, \g_v$ implies 
$\dim L_v=\frac{1}{2}\dim S_p$. Let us check that $L_v$ is isotropic. 
For $s,t=\g_vt'\in L_v={\rm im}\, \g_v$, we have 
\[ h(s,t)=h(s,\g_vt')=\t h(\g_vs,t')=0,\] 
since $s\in L_v=\ker \g_v$. \qed

\noindent 
The lemma shows that ${\rm rk}\, \g_v =
\frac{1}{2}\dim S_p$ for any non-zero null vector. 
Now we consider the bilinear form $\beta =h(\g_v \cdot ,\cdot )$ 
on $S_p$; ${\rm rk}\, \beta = {\rm rk}\, \g_v =\frac{1}{2}\dim S_p$. 
Under the assumption $\g_v S_0 \subset S_0^\perp$, the matrix of $\beta$ 
with respect to a basis adapted to a direct decomposition
$S_p=S_0\oplus S_1$ is of the form
\[ \left( \begin{array}{cc} 0&A\\
\s \t A^t&B
\end{array}\right) \] 
(Notice that the symmetry of $\beta$ is $\s \t$.) Therefore, 
\[ \frac{1}{2}\dim S_p = {\rm rk}\, \beta \le  {\rm rk}\, A + 
 {\rm rk}\, (\s \t A^t,B) \le 2 \dim S_1=2(\dim S_p-\dim S_0),\]
which implies $\dim S_0 \le \frac{3}{4}\dim S_p$. So 
$\dim S_0>\frac{3}{4}\dim S_p$ implies 
$\nexists v\in T_pM\setminus \{ 0\}: \g_v S_0\subset S_0^\perp$.
This shows that $\Pi : S_0\ot S_0 \ra T_pM$ is surjective in case (i) and 
(ii). 

The proof of (iii) uses the fact that in that case 
there exist a unique (up to a constant factor) 
admissible parallel nondegenerate bilinear form $h$ invariant under
$J_1$, $J_2$ and $J_3$, see \cite{AC}. The form is of type $\t=+1$. 
Using this form we obtain 
for two solutions $s, t$ of \re{imagKSEqu} that 
$Y=\o =[s,t]_1$ is a conformal Killing vector field, which satisfies 
\re{confKvEqu} with 
\[ \tilde{\o} = -h(\tilde{s},t)- h(s,\tilde{t})= 
-\l (h(Is,t)+h(s,It))=0,\] 
as follows from \re{tildeEqu}. Therefore $Y$ is a Killing vector field. 
The rest of the proof is similar to that of (i) and (ii).
   
To prove (iv) we first remark that the assumptions on the dimension and 
signature ensure the existence of an admissible parallel nondegenerate 
bilinear form $h$ of type $\t=+1$. Then we consider the subspaces 
$S_0(\l ), S_0(-\l )\subset S_p$
spanned by the values at $p$ of Killing spinors with Killing numbers
$\l$ and $-\l$, respectively. In virtue of Corollary \ref{mainCor},
$[s,t]_1$ is a Killing vector field if $s$, $t$ are Killing spinors
with Killing numbers $\l$, $-\l$, respectively.  Therefore, it suffices
to show that $[S_0(\l ),S_0(-\l )]=T_pM$. If this condition were not
fulfilled, there would exist $0\neq v\in T_pM$ such that
$\g_v : S_0(\l) \ra S_0(-\l)^\perp$.  The assumption $\dim S_0(-\l)^\perp = 
N-k_-<k_+=\dim S_0(\l)$ implies that $L_v=\ker \g_v \neq 0$. Then $g$ is 
indefinite,  
$v$ is a null vector and  $L_v={\rm im}\, \g_v$ is maximally isotropic,
by Lemma \ref{Lemma}. In particular, ${\rm rk}\, \g_v =N/2$. We can consider
$\b = h(\g_v \cdot ,\cdot )$ as a linear map
$S_p \ra S_p^*$. {}From the matrix representation of $\b$ with respect
to bases adapted to decompositions 
$S_p=S_0(\l ) \oplus S_1$ and $S_p^*\cong S_0(-\l)^*\oplus 
S_1'$ we see that 
\[ \frac{1}{2}N={\rm rk}\, \b \le \min (k_-,N-k_+)+N-k_-=2N-k_+-k_-, \]
and, hence, $k_++k_-\le \frac{3}{2}N$, which contradicts the assumption 
$k_++k_-> \frac{3}{2}N$. This proves $[S_0(\l ),S_0(-\l )]=T_pM$.
\qed

Now we study the case where the bilinear form $h$ has type $\t=+1$
and the Killing spinors have the same Killing number. 
\bp Let $h$  be a nondegenerate parallel bilinear form
of symmetry $\s$ and type $\t =+1$ on the spinor bundle $S$ 
of a pseudo-Riemannian spin manifold $(M,g)$ and denote by 
$S(\l )\subset \G (S)$ the
vector space of Killing spinors with a given Killing number $\l \in \bR
\setminus\{ 0\}$.
Then the image $[S(\l ),S(\l )]_1\subset \G (TM)$ consists of 
Killing vector fields if and only if $S_0:=S(\l )|_p\subset S_p$ 
is totally isotropic for all $p\in M$ with respect to $h$.  If $S_0$ is 
maximally isotropic at a point $p$ then $[S(\l ),S(\l )]_1\neq 0$,
hence, $(M,g)$ admits a Killing vector field, which does not vanish at $p$.
\ep 

\pf By Corollary \ref{mainCor}, the bracket $\o = [s,t]_1$ of $s, t\in S(\l )$
is a conformal Killing vector field with 
\[ \tilde{\o} = -2\l h(s,t).\] 
This shows that $\o$ is a Killing vector field if and only if $h(s,t)=0$. 
Assume now that $S_0=S_0^\perp$ is 
maximally isotropic. By \re{fundEqu}, 
$[S_0,S_0]_1=0$ is equivalent to $\g_vS_0\subset S_0$ 
for all $v\in T_pM$, which is impossible since $S_p$ is an irreducible
module of the Clifford algebra $C\!\ell (T_pM)$. \qed  

\noindent 
{\bf Remark:} One can check that $[S_0,S_0]_1$ is one-dimensional
for any maximally isotropic subspace $S_0$ of the spinor module
$S_{2,3}=\bR^4$ of $\Spin(2,3)$. For the spinor module $S_{4,5}$ of
$\Spin(4,5)$ one can construct a maximally isotropic subspace
$S_0$ such that $\dim [S_0,S_0]_1=4$. These examples show that in general
a vector space of Killing spinors spanning a maximally isotropic
subspace of $S_p$ for all $p$ is not sufficient to
produce a transitive Lie algebra of Killing fields.      
\section{A multiplicative invariant}
Let $M$ be a pseudo-Riemannian spin manifold with real spinor bundle $S$
of rank $N$ and denote by $S(\l )=S(M,\l )$ the vector space of Killing spinors 
with Killing number $\l \in \bR$. Then we put $k:= \dim S(\l)$ and 
\[ \k (M,\l ) := \frac{k}{N},\quad \k (M):=\k (M,0).\] 
Notice that $\k (M)=1$ if and only if $M$ is flat and that
$\k (M,\l )= \frac{\dim \bS (\l)}{{\rm rk}\, \bS}$, where $\bS$ is the 
complex spinor bundle and $\bS (\l )=\bS (M,\l )$ the vector space
of complex Killing spinors with Killing number $\l$. This follows from the fact that the complex spinor module $\mathbb{S}_ {p,q}$ of $C\!\ell_{p,q}$
is either the complexification of the real spinor module $S_{p,q}$
or coincides with $S_{p,q}$ endowed with a Pin$(p,q)$-invariant 
complex structure, see \cite{ACDV} Table 1.
As a consequence, we have ${\rm rk}\, \bS = N$ or $N/2$, respectively. 

\bl \label{L} 
Let $V=V_1+V_2$ be an orthogonal decomposition of a complex Euclidian
vector space of dimension $n$ into subspaces of dimension $n_1$, $n_2$
respectively.  
\begin{enumerate}
\item[(i)] If $n_1$ or $n_2$ is even, then the Clifford algebra
$C\!\ell (V)\cong C\!\ell (V_1)\ot C\!\ell (V_2)$ and 
the tensor product $S(V)=S(V_1)\ot S(V_2)$ of irreducible
$C\!\ell (V_1)$-,  $C\!\ell (V_2)$-modules $S(V_1)$ and $S(V_2)$, respectively, 
is an irreducible $C\!\ell (V)$-module. 
\item[(ii)] If $n_1$ and $n_2$ are odd, then 
$C\!\ell (V)\not\cong C\!\ell (V_1)\ot C\!\ell (V_2)$ but 
$C\!\ell (V)$ is isomorphic to the $\bZ/2\bZ$-graded tensor product  
$C\!\ell (V)\cong C\!\ell (V_1)\hat{\ot} C\!\ell (V_2)$.
In this case the spinor module of $C\!\ell^0 (V)$ is obtained
as the even part $(\S \hat{\ot} \S')_0=\S_0\hat{\ot}\S_0' +
\S_1\hat{\ot}\S_1'$ of the $\bZ/2\bZ$-graded tensor product 
$\S\ot \S'$ of irreducible 
$\bZ/2\bZ$-graded  $C\!\ell (V_1)$-,  
$C\!\ell (V_2)$-modules $\S , \S'$, respectively. The $C\!\ell^0 (V)$-module 
$S(V)=(\S \hat{\ot} \S')_0$ is a sum of non-equivalent irreducible
semi-spinor submodules $S^\pm(V)$, which are the $\pm i$-eigenspaces
of a central element $\xi \in C\!\ell^1(V_1)\hat{\ot}C\!\ell^1(V_2)$  
of $C\!\ell^0 (V)$. 
\end{enumerate}
\el 
\bc \label{prodCor} Under the assumptions of Lemma \ref{L} the following is true.
\begin{enumerate}
\item[(i)] If $n_1$ or $n_2$ is even, then as a  
$\Spin(V_1)\times \Spin(V_2)$-module the spinor module $S(V)$  
of $\Spin(V)$ is isomorphic to the tensor product
$S(V)\cong S(V_1)\ot S(V_2)$.
\item[(ii)] If $n_1$ and $n_2$ are odd, then as a  
$\Spin(V_1)\times \Spin(V_2)$-module, $S(V)\cong 2S(V_1)\ot S(V_2)$.
\end{enumerate}
\ec

\bc \label{kappaCor} 
Let $M=M_1 \times M_2$ be the product of two pseudo-Riemannian
spin manifolds. Then 
$\k (M)=\k (M_1) \k (M_2)$. In particular,
$\k (M)=\k (M_1)$ if and only if $M_2$ is flat. 
\ec 

\pf Since $\k (M) = \frac{\dim \bS (M,0)}{{\rm rk}\, \bS}$, the 
statement of Corollary \ref{kappaCor} 
is obtained from Corollary  \ref{prodCor},  
using that parallel spinors correspond to invariants of the 
holonomy group under the spinor representation and that 
the holonomy group of $M$ is the product of
the holonomy groups of the factors $M_1, M_2$. In fact,
$\bS (M,0)\cong \bS (M_1,0)\ot \bS (M_2,0)$ if $n_1$ and $n_2$ are
even and $\bS (M,0)\cong 2\bS (M_1,0)\ot \bS (M_2,0)$ if $n_1$ and $n_2$ 
are odd.  
\qed 

\noindent
{\bf Remark:} The invariant $\k (M,\l )$ for $\l\neq 0$ is not multiplicative.
For instance, $\k (S^2,\frac{1}{2})=1$ but 
$\k (S^2\times S^2,\frac{1}{2})=0$.  
\section{Manifolds with many parallel spinors}
\bt \label{parallelThm} Let $(M,g)$ be a pseudo-Riemannian spin manifold. 
\begin{enumerate} 
\item[(i)] If $\k (M)>\frac{3}{4}$,
then $(M,g)$ is flat. 
\item[(ii)] If the metric $g$ is definite and 
$\k (M) > \frac{1}{4}$, then  $(M,g)$ is flat. 
A complete simply connected 
Riemannian spin manifold $(M,g)$ with $\k (M) = \frac{1}{4}$ 
is 
the product of a flat manifold and a manifold with holonomy group $\SU (2)$. 
\end{enumerate}
\et  

\pf (i) follows from Theorem \ref{2ndThm} (i) and (v), since the  conformal
Killing vector fields $[s,t]_1$ are parallel if $s, t$ are parallel spinors,
see Corollary \ref{mainCor}.\\
Next we prove (ii). It follows from 
Wang's classification of parallel spinors on manifolds 
with connected irreducible holonomy group \cite{W} 
that a locally irreducible Riemannian manifold $(M,g)$ has 
$\k (M) \le \frac{1}{4}$ and $\k (M) =\frac{1}{4}$
implies that $M$ has holonomy algebra $\gh =\su (2)$. 
Applying the (local) de Rham decomposition and 
Corollary \ref{kappaCor}, we conclude that 
a Riemannian manifold with $\k (M)>\frac{1}{4}$ is flat and that 
a complete simply connected Riemannian manifold with
$\k (M) = \frac{1}{4}$ is 
the product $M=M_0\times M_1$ of a flat manifold $M_0$ and an irreducible 
manifold $M_1$ with $\k (M_1) = \frac{1}{4}$ and holonomy group 
$\SU (2)$.  
\qed 

\bt \label{LorCone} Let $(\hat{M},\hat{g})$ be the Lorentzian cone over a 
pseudo-Riemannian
manifold $(M,g)$ with either negative definite metric 
or of signature $(+,\ldots ,+,-)$. If $\k (\hat{M})>\frac{1}{2}$, then
$\hat{M}$ is flat and $M$ has constant curvature $1$. 
\et

\pf  If $(\hat{M},\hat{g})$ is not flat, we can decompose it  
locally as a product of
indecomposable pseudo-Riemannian manifolds. By Corollary 
\ref{kappaCor}, there exists an indecomposable factor $M_1$
of dimension $>1$ with $\k (M_1) > \frac{1}{2}$. It cannot
be Riemannian, by the  previous theorem. Hence it is a
Lorentzian indecomposable manifold. By \cite{ACGL} Theorem 4.1,
$M_1=\hat{N_1}$ is (locally) a cone over a pseudo-Riemannian manifold $N_1$. 
Moreover, by \cite{ACGL} Theorem 9.1, the local holonomy algebra $\hat{\gh}$ 
of $M_1$ contains the subalgebra $\mathfrak{e} := p\wedge E$, where 
$T_xM_1=V=\bR p+\bR q+E$,  
$p, q$ are isotropic vectors with $\hat{g}(p,q)=1$ and 
$E$ is the positive definite orthogonal
complement of ${\rm span}\{ p,q\}$. The Clifford algebra 
has the decomposition $C\!\ell (V)=C\!\ell_{1,1}\ot C\!\ell (E)$.
The Clifford algebra $C\!\ell_{1,1}$ is the full matrix algebra
of real $2\times 2$ matrices and is generated by 
\[ \g_p = \sqrt{2}\left( \begin{array}{cc}0&1\\
0&0
\end{array}
\right),\quad \g_q=-\sqrt{2}\left( \begin{array}{cc}0&0\\
1&0
\end{array}
\right),\]
with respect to the standard basis $(e_1,e_2)$ of $\bR^2$.  
If $S_E$ is an irreducible $C\!\ell (E)$-module then
$S_V=\bR^2\ot S_E$ is an irreducible $C\!\ell (V)$-module. 
Under the isomorphism $C\!\ell (V)\cong C\!\ell_{1,1}\ot C\!\ell (E)$ a vector 
$v=f\oplus e\in \bR^{1,1}\oplus E$ is mapped to $f\ot 1 + \nu \ot e$,
where $\nu = \frac{1}{2}(p+q)(p-q)$ is the volume element in $C\!\ell_{1,1}$,
which satisfies $\nu^2=1$, $\nu e_1=e_1$, $\nu e_2=-e_2$. 

\bl \label{invL} The space of $\mathfrak{e}$-invariant spinors is given by
\[ S_V^\mathfrak{e} = e_1\ot S_E\subset S_V=\bR^2\ot S_E.\]  
\el 

\pf A spinor $s=e_1\ot s_1 + e_2\ot s_2\in S_V$
is invariant under $\mathfrak{e} \subset \hat{\gh}$ if and only if
\[ 0=\g_{p\wedge e} s= -\sqrt{2} e_1\ot \g_es_2,\]
for all $e\in E$, which is equivalent to $s_2=0$.
\qed  

\noindent 
The lemma shows that $\dim S_V^{\hat{\gh}} \le \dim S_V^\mathfrak{e} = 
\frac{1}{2}\dim S_V$ and, hence, $\k (M_1) \le \frac{1}{2}$,
which contradicts the assumption. 
\qed 
\section{Cones $\hat{M}$ over pseudo-Riemannian manifolds $M$ and
relation between Killing spinors on $M$ and parallel spinors on $\hat{M}$}
\bd \label{coneDef} 
Let $(M,g)$ be a pseudo-Riemannian manifold of signature $(p,q)$. 
The manifold $\hat{M}=\bR^+ \times M$ endowed with the  pseudo-Riemannian
metric $\hat{g}= dr^2+r^2g$ of signature $(p+1,q)$ is called the {\cmssl cone}
over $(M,g)$.  
\ed 
Recall that a {\cmssl spin structure} (in the strong sense)  
on $(M,g)$ is 
a $\Spin_0(p,q)$-equivariant two-fold covering 
$P_{\Spin_0(p,q)}(M) \ra P_{SO_0(p,q)}(M)$ of the 
principal bundle of strongly oriented orthonormal frames. Let us denote
by $P_{\Spin_0(p+1,q)}(M) \supset P_{\Spin_0(p,q)}(M)$ and  $P_{SO_0(p+1,q)}(M) 
\supset P_{SO_0(p,q)}(M)$ the $\Spin_0(p+1,q)$- and  $SO_0(p+1,q)$-principal
bundles obtained by enlarging the structure groups. Then 
$P_{\Spin_0(p,q)}(M) \ra P_{SO_0(p,q)}(M)$ extends naturally to a 
$\Spin_0(p+1,q)$-equivariant two-fold covering
\[ \Theta : P_{\Spin_0(p+1,q)}(M) \ra P_{SO_0(p+1,q)}(M).\] 
Using the isometric inclusion $M\cong \{ 1\}\times M \subset 
\hat{M}= \bR^+\times M$, we can identify $P_{SO_0(p+1,q)}(M)$ 
with the restriction $P_{SO_0(p+1,q)}(\hat{M})|_M$ of the 
bundle of strongly oriented orthonormal frames of $\hat{M}$. 
In particular, the frame $(e_1,\ldots ,e_n)\in P_{SO_0(p,q)}(M)\subset 
P_{SO_0(p+1,q)}(M)$
is mapped to the frame $(\partial_r,e_1,\ldots ,e_n)
\in P_{SO_0(p+1,q)}(\hat{M})$ under
this identification. 
Similarly, we identify $P_{SO_0(p+1,q)}(\hat{M})$ with the 
pullback of $P_{SO_0(p+1,q)}(M)$ via the projection $\pi : \hat{M} \ra M$. 
Under this identification $(\partial_r,e_1,\ldots ,e_n)
\in P_{SO_0(p+1,q)}(\hat{M})_{(r,x)}$ is mapped to  
$(re_1,\ldots ,re_n)\in P_{SO_0(p,q)}(M)_x\subset P_{SO_0(p+1,q)}(M)_x$ 
for all $x\in M$.
Then 
\[P_{\Spin_0(p+1,q)}(\hat{M}):=\pi^*P_{\Spin_0(p+1,q)}(M) 
\ra \pi^*P_{SO_0(p+1,q)}(M)=
P_{SO_0(p+1,q)}(\hat{M})\] 
defines a spin structure on $\hat{M}$.  

\bl  \label{realL} Let $(\hat{M},\hat{g})$ be the cone over a 
pseudo-Riemannian spin manifold $(M,g)$ of signature $(p,q)$. 
\begin{enumerate}
\item[(i)] If $s=p-q\equiv 0,2,4,5$ or $6\pmod{8}$, then the 
spinor bundle  $\hat{S}$ of $\hat{M}$ is related to the 
spinor bundle $S$ of $M$ by a canonical isomorphism  
\[ \hat{S}|_M \cong S.\] 
\item[(ii)] If $s=p-q\equiv 1, 3$ or $7\pmod{8}$, then the 
semi-spinor bundles  $\hat{S}^\pm$ of $\hat{M}$ are related to the 
spinor bundle of $M$ by canonical isomorphisms   
\[ \hat{S}^\pm|_M \cong S.\] 
\item[(iii)] If $n=\dim M =p+q$ is even, then the complex spinor bundles
$\bS$, $\hat{\bS}$ of $M$ and $\hat{M}$, respectively, 
are related by a canonical
isomorphism
\[ \hat{\bS}|_M \cong \bS.\]
\item[(iv)] If $n$ is odd, then the complex semi-spinor bundles $\hat{\bS}^\pm$
of $\hat{M}$ are  related to the spinor bundle $\bS$ of $M$ by  canonical
isomorphisms
\[ \hat{\bS}^\pm|_M \cong \bS.\]
\end{enumerate}
\el 

\pf Let $(e_0,\ldots ,e_n)$ be an orthonormal basis of
$\bR^{p+1,q}$. Recall that by definition  
\[ \Spin(p,q) \subset \Spin (p+1,q) \subset 
C\!\ell^0_{p+1,q}=\langle e_ie_j| i,j=0,\ldots ,n\rangle.\] 
The even part $C\!\ell^0_{p+1,q}$ of the
Clifford algebra $C\!\ell_{p+1,q}$ is mapped isomorphically
onto $C\!\ell_{p,q}$ by
\begin{eqnarray*}  e_ie_j &\mapsto& e_ie_j,\\
e_ie_0&\mapsto& e_i,\quad i,j= 1,\ldots ,n.
\end{eqnarray*}
Using this isomorphism $C\!\ell^0_{p+1,q}\cong C\!\ell_{p,q}$ 
the spinor module $S_{p,q}$ of $\Spin(p,q)$ can be 
extended to an irreducible $\Spin(p+1,q)$-module. In fact, $S_{p,q}$
is the restriction of an irreducible $C\!\ell_{p,q}$-module 
to $\Spin (p,q) \subset C\!\ell_{p,q}$. Restricting this 
$C\!\ell_{p,q}$-module to $\Spin (p+1,q)\subset  C\!\ell^0_{p+1,q}\cong 
C\!\ell_{p,q}$ gives the desired $\Spin(p+1,q)$-module. The 
$\Spin(p+1,q)$-module $S_{p,q}$ is equivalent to the spinor module
$S_{p+1,q}$ if the spinor module
$S_{p+1,q}$ is irreducible, which is the case if 
$s=p-q\equiv 0,2,4,5$ or $6\pmod{8}$, see \cite{AC} Prop.\ 1.3. 
Otherwise it is equivalent to one of the semi-spinor modules
$S^\pm_{p+1,q}$. The semi-spinor modules $S^+_{p+1,q}$ and $S^-_{p+1,q}$ are always
equivalent as $\Spin(p,q)$-modules (and even as $\Spin(p+1,q)$-modules
if $s\equiv 1\pmod{8}$). This implies (i) and (ii). The proof of (iii) 
and (iv) is similar. 
\qed 

Notice that for $s=p-q\equiv 5\pmod{8}$ the spinor module $S_{p+1,q}$ 
is irreducible and admits a $\Spin(p+1,q)$-invariant complex structure, 
see \cite{AC} Prop.\ 1.3. 
Its complexification is isomorphic to the complex spinor module
of $\Spin(p+1,q)$ (see \cite{ACDV} Table 1), 
which is a sum of two semi-spinor modules. 

For $\S \in \{ S, \bS , \hat{S}, \hat{\bS}, \hat{S}^\pm, \hat{\bS}^\pm\}$, 
let us denote by $\S  (\l )$ the vector space of  
Killing spinors $s\in \G (\S )$ with Killing number $\l \in \bR$. 

Notice that if $\l\neq 0$ one can always normalise the metric such that
$\l=\pm \frac{1}{2}$ (as for a space of constant curvature 1). Now let $\S = S$ or $\bS$. 
Multiplication
by the volume element $\nu =e_1\cdots e_n\in C\!\ell (TM)$ 
maps $\S (\l )$ to $\S ((-1)^{n+1}\l )$. 
In particular, it defines  
isomorphisms  $\S (\l )\cong \S (-\l )$, 
if $n$ is even. For
odd dimensional manifolds, however, the  vector spaces 
$\S (\l )$ and $\S (-\l )$ have in general different dimensions.  

Using Lemma \ref{realL},  
the following theorem can be proven as for Riemannian manifolds, see 
B\"ar \cite{B}. 
\bt \label{BaerThm} Let $(\hat{M},\hat{g})$ be the cone over a 
pseudo-Riemannian spin manifold $(M,g)$ of signature $(p,q)$. 
\begin{enumerate}
\item[(i)] 
The restriction $\G (\hat{S})\ni s\mapsto s|_M\in \G (S)$ 
defines isomorphisms
\[ \hat{S} (0)\ra  S \left(\frac{1}{2}\right)\cong  S 
\left(-\frac{1}{2}\right),\]
if $s=p-q\equiv 0,2,4,5$ or $6\pmod{8}$ and  
\[ \hat{S}^\pm (0)\ra S \left(\pm \frac{\e }{2}\right),\]
for some $\e \in \{ 1,-1\}$, if $s=p-q\equiv 1,3$ or $7\pmod{8}$. 
\item[(ii)] The restriction $\G (\hat{\bS})\ni s\mapsto s|_M\in \G (\bS)$ 
defines isomorphisms
\[ \hat{\bS} (0)\ra  \bS \left(\frac{1}{2}\right)\cong  \bS 
\left(-\frac{1}{2}\right),\]
if $n=\dim M$ is even and 
\[ \hat{\bS}^\pm (0)\ra \bS \left(\pm \frac{\e }{2}\right),\]
for some $\e \in \{ 1,-1\}$, if $n$ is odd.  
\end{enumerate}
\et 

\section{Riemannian manifolds with many Killing spinors}
\bt \label{mixThm} 
Let $(M,g)$ be a simply connected Riemannian spin manifold.  
\begin{enumerate}
\item[(i)] 
Assume that one of the following conditions is satisfied:
\begin{enumerate}
\item[a)] $(M,g)$ is complete and not of constant
curvature $1$.
\item[b)] The holonomy algebra of $M$ is different 
from $\mathfrak{so} (n)$, where 
$n=\dim M$.
\end{enumerate}
Then the holonomy algebra $\hat{\mathfrak{h}}$ of the cone 
$(\hat{M},\hat{g})$ is irreducible. 
\item[(ii)] 
If $\hat{\mathfrak{h}}$ is irreducible, then  
$(\hat{M},\hat{g})$ admits 
a parallel spinor if and only if $\hat{\mathfrak{h}}$ 
belongs to the following list: 
$\mathfrak{su}(m)$ $(m\ge 3, k=2)$, $\mathfrak{sp}(m)$ $(m\ge 2, k=m+1)$, 
$\mathfrak{spin} (7)$ $(k=1)$ 
or $\mathfrak{g}_2$ $(k=1)$, where $k$ in brackets indicates 
the number of linearly independent 
parallel complex spinors.  The projection of the space of 
parallel complex spinors onto the space of parallel complex semi-spinors 
is zero for one of the two semi-spinor bundles, 
unless $n+1\not\equiv 0\pmod{4}$. 
\end{enumerate}
\et

\pf  The irreducibility of the holonomy algebra follows from
Gallot's theorem \cite{G} under the assumption a) and from 
\cite{ACGL} Theorem 4.1 under the assumption b). The remaining
statements follow from 
Wang's classification
of parallel spinors on 
manifolds with connected irreducible holonomy group \cite{W}
and the observation that there is no cone with holonomy group 
$\SU (2) = \Sp (1)$. 
\qed 

\bt \label{RiemThm} 
Let $(M,g)$ be a Riemannian spin manifold 
which is not of constant 
positive curvature $4\l^2$. 
\begin{enumerate}
\item[(i)] 
Then $\k (M,\l ) \le  \frac{3}{8}$. 
\item[(ii)] Assume that for every $p\in M$ we have 
$\k (U,\l ) =  \frac{3}{8}$ for every sufficiently 
small open neighborhood $U\subset M$ of $p$. 
Then either $(M,g)$ is locally isometric to seven-dimensional 
$3$-Sasakian manifold or there exists a dense open subset
$M'\subset M$ such that every point of $M'$ has a neighborhood isometric to 
a Riemannian manifold of the form 
\[
(I \times M_1 \times M_2, ds^2 + \cos^2(s) g_1 + \sin^2(s) g_2),
\]
where $(M_1,g_1)$ is of constant curvature $1$ or of dimension $\le 1$ 
and $(M_2,g_2)$ is a seven-dimensional $3$-Sasakian manifold.
\end{enumerate}
\et

\pf  We use the correspondence between Killing spinors on $M$ and
parallel spinors on the cone $\hat{M}$. If $\k (M ,\l ) >\frac{3}{8}$,
then, according to Theorem \ref{BaerThm}, $\k (\hat{M}) > 
\frac{1}{2} \times \frac{3}{8} =\frac{3}{16}$. This is impossible
if the holonomy algebra $\hat{h}$ of $\hat{M}$ is irreducible, due to
Theorem \ref{mixThm} (ii). The maximal value
$\k (\hat{M}) = \frac{3}{16}$ is, in fact, attained for the holonomy algebra
$\hat{\gh}= \gsp (2)$ (since there is no cone with holonomy 
$\su (2)$). The cone $\hat{M}$ has local holonomy $\gsp (2)$ if and only if
the seven-dimensional manifold $M$ is locally 3-Sasakian. 
This proves (i) and (ii) if $\hat{\gh}$ is irreducible. 
In the reducible case, the claims (i) and (ii) now follow from 
\cite{ACGL} Theorem 4.1 using Corollary \ref{kappaCor}.   
\qed  

Recall \cite{B} that the holonomy algebra $\hat{\mathfrak{h}}$ of the 
cone $(\hat{M},\hat{g})$ over a simply connected Riemannian 
manifold $(M,g)$ belongs to the list of irreducible linear Lie algebras 
described in Theorem \ref{mixThm} (ii) 
if and only if $(M,g)$ is Einstein-Sasaki,  $3$-Sasakian, strictly 
nearly parallel ${\rm G}_2$ 
or strictly nearly K\"ahler, respectively. We will call  these  
geometric structures on $(M,g)$ {\cmssl B\"ar geometries}. 
\bt \label{ClassThm} Let $(M,g)$ be an $n$-dimensional Riemannian spin manifold
which admits a nontrivial Killing spinor with Killing constant $\l\in \bR$.
\begin{enumerate} 
\item[(i)] If $\l=0$, then $(M,g)$ is locally a product 
$M=M_0\times M_1 \times \cdots \times 
M_r$ of a flat Riemannian manifold $M_0$ with an arbitrary number  
of Riemannian manifolds $M_i$ with irreducible holonomy group from the 
following list: $\SU (m)$, $\Sp (m)$, 
$\Spin (7)$ or ${\rm G}_2$. 
\item[(ii)] If $\l \neq 0$, then  $(M,g)$ has holonomy $\gh=\so (n)$. 
Moreover, if the cone $\hat{M}$ is locally irreducible, then
$(M,g)$ carries locally 
one of the B\"ar geometries and if $\hat{M}$ is locally reducible, 
then, on a dense open subset, $(M,g)$   
can be locally 
represented in the 
form
\begin{equation}\label{decompEqu} 
 M= I \times M_1 \times M_2,\quad g=ds^2 + \cos^2(s) g_1 + \sin^2(s) g_2,
\end{equation}
where $I \subset (0,\frac{\pi}{2})$ is an intervall and 
$(M_1,g_1)$ and $(M_2,g_2)$   
are Riemannian manifolds which either  
admit a nontrivial Killing spinor with 
Killing constant $\pm \l$ or which are of dimension $\le 1$.   
\end{enumerate}
\et 

\noindent
{\bf Remark:}
Notice that the Theorem \ref{ClassThm} (ii) 
gives an inductive decomposition of a manifold 
with a nontrivial Killing spinor with $\l\neq 0$ in terms of an 
arbitrary number of manifolds $(M_i,g_i)$, which carry each one of 
the B\"ar geometries or are of dimension
$\le 1$. We remark also that a manifold $(M,g)$ of constant curvature $1$
can be locally decomposed as  \re{decompEqu} with manifolds $(M_1,g_1)$ and 
$(M_2,g_2)$ which are either  
of constant curvature $1$ or of dimension $\le 1$.  

\pf  (i) is an immediate consequence of Wang's classification of irreducible
connected holonomy groups preserving a non-trivial spinor \cite{W}.

In the case (ii), it follows from B\"ar's classification \cite{B} that 
the cone $(\hat{M},\hat{g})$ over $(M,g)$ is locally irreducible if and only 
if $(M,g)$ carries locally one of the B\"ar geometries. 
We check, in this case, that $M$ is locally irreducible. 
We consider the modified covariant derivative $\tilde{\n}_X := \n_X - \l \g_X$, 
$X\in T_xM$,  on the complex spinor bundle over $M$, where $\n$ stands for the
Levi-Civita connection. 
Assume that $M=M_1\times M_2$ is a Riemannian product. Then we compute 
the curvature of $\tilde{\n}$ at $x\in M$:
\[ \tilde{R}(X_1,X_2)_x = \l [\g_{X_1},\g_{X_2}] = 2\l \g_{X_1}\g_{X_2}, \] 
for $X_i\in TM$
tangent to $M_i$ and such that $\n X_i|_x=0$, $i=1,2$. This 
implies that the local holonomy $\tilde{\gh}$ algebra of $\tilde{\n}$ contains 
$\spin (n)$, because the holonomy algebra at $x$ contains all 
curvature operators $\tilde{R}(X_1,X_2)_x$ and the Clifford products 
$X_1 X_2$ generate $\spin (n)$ (as a Lie algebra). 
Since, by \cite{B}, $\tilde{\gh}$ can be identified 
with the local holonomy algebra 
$\hat{\gh}$ of the Levi-Civita 
connection of the cone
$\hat{M}$,  we can conclude that $\hat{\gh}$ 
contains the subalgebra $\so (n)\subset \so (n+1)$. 
One can easily check that this is not possible
for $\hat{\gh}$ belonging to the list of irreducible holonomy algebras
of Riemannian cones admitting a parallel spinor, see   
Theorem \ref{mixThm} (ii). This shows that $(M,g)$ is locally irreducible
if $(\hat{M},\hat{g})$ is locally irreducible. 
In particular, $\gh$ belongs to Berger's list of irreducible holonomy algebras, 
excluding the Ricci-flat holonomies (but 
so far including the holonomies of irreducible symmetric spaces). 
In dimension $n=7$ this already implies
$\gh =\so (7)$. In dimension $n=6$ this implies $\gh =\so (6)$, using that 
a strict nearly K\"ahler manifold cannot be K\"ahler. In the remaining cases
$(M,g)$ is locally Sasaki-Einstein (or even 3-Sasakian). {}The curvature tensor
of such a manifold satisfies
\[ R(\xi ,X)Y= \xi g(X,Y) - Xg(\xi ,Y)\] 
for all vector fields $X, Y$ on $M$, where $\xi$ is the Sasaki vector field. 
This identity immediately implies that $\gh = \so (n)$, since 
$\gh$ contains all curvature operators and their brackets.    

If the cone $(\hat{M},\hat{g})$
is locally reducible, then it  
follows from \cite{ACGL} Theorem 4.1 that $\gh=\so (n)$ and that, on a dense
open subset of $M$, $(M,g)$ is locally isometric to \re{decompEqu}.  
\qed

\section{Pseudo-Riemannian manifolds with Lorentzian cone,
which admit many Killing spinors} 
\bt \label{Lor/ndThm} Let $(M,g)$ be spin with either a 
 negative definite metric or a metric of Lorentzian signature 
$(+,\ldots ,+,-)$. If $(M,g)$ is not of positive 
constant curvature $4\l^2$, then $\k(M,\l ) \le \frac{1}{2}$. 
\et 

\pf 
The spinor module $S_V$ of $\Spin(V)$, $V=T_x\hat{M}$, is either irreducible 
or it  splits as 
\[ S_V=S_V^+ \oplus S_V^-,\quad 
S_V^\pm = e_1\ot S_E^\pm + e_2\ot S_E^\mp, 
\]
where $S_E^\pm$ are the semi-spinor modules of $\Spin(E)$ and we use the 
notation of Lemma \ref{invL}. As in the proof of Theorem \ref{LorCone}, we
can assume that the cone $(\hat{M},\hat{g})$ is indecomposable. 
If the cone is Riemannian, we have $\k (M,\l) \le \frac{3}{8}$, by
Theorem \ref{RiemThm}. Therefore, we can assume that it is Lorentzian. 
In that case  the holonomy algebra contains $\mathfrak{e}=p\wedge E$, by 
\cite{ACGL} Theorem 9.1. 
Then $(S_V)^\mathfrak{e}= e_1\ot S_E$ if $S_V$ is irreducible and 
$(S_V^\pm)^\mathfrak{e}= e_1\ot S_E^\pm$ otherwise. This shows that 
$\dim \hat{S}(0) =\frac{1}{2}{\rm rk} \hat{S}$, which implies  
$\dim S (\l ) \le \frac{1}{2}N$, by Theorem \ref{BaerThm}.  
\qed 

Remark that a pseudo-Riemannian manifold $(M,g)$ of dimension $n$ which admits 
a Killing spinor with (real) Killing number 
$\l\in \bR \setminus \{ 0\}$ has positive scalar curvature  $s=4n(n-1)\l^2$. 
If $g$ is negative definite of scalar curvature 
$s>0$, then the Riemannian
metric $-g$ has negative scalar curvature $-s$. This allows to 
treat also Riemannian manifolds with negative scalar curvature. 
 
\end{document}